# POLARIZATIONS ON ABELIAN VARIETIES AND SELF-DUAL $\ell$-ADIC REPRESENTATIONS OF INERTIA GROUPS

A. SILVERBERG AND YU. G. ZARHIN

1. INTRODUCTION

In this paper we obtain isogeny classes of abelian varieties all of whose polarizations have degree divisible by a given prime (for all odd primes $\ell$ and for all dimensions $d \geq \ell - 1$). In particular, these abelian varieties do not admit a principal polarization. Everett Howe ([8]; see also [9]) obtained examples of isogeny classes of abelian varieties with no principal polarizations. We believe our results give the first examples for which all the polarizations of the abelian varieties in the isogeny classes have degree divisible by a given prime. Howe's examples are of ordinary abelian varieties over finite fields, while ours are supersingular abelian varieties over infinite fields (in every positive characteristic).

These results rely on a study of $\ell$-adic representations of finite groups, especially inertia groups (see §4–§6). In recent years, $\ell$-adic Galois representations associated to abelian varieties have played an important role in answering number-theoretic questions. One motivation for the representation-theoretic results in this paper was to better understand the $\ell$-adic symplectic representations of inertia groups which arise from abelian varieties.

We will say that a group is "of inertia type" if it is a finite group which has a normal Sylow-$p$-subgroup with cyclic quotient, for some prime $p$ ($\neq \ell$). These groups are exactly the inertia groups of finite Galois extensions of discrete valuation fields of residue characteristic $p$ (see Corollaire 4 and the remark which follows in Chapitre IV, §2 of [16]).

In [22] we studied certain groups of inertia type that arose from work of Grothendieck and Serre. These groups are the Galois groups of the smallest extensions over which abelian varieties over local fields acquire semistable reduction. We showed that these groups embed in a product of the $\mathbf{Q}_\ell$-points of a symplectic group and the $\mathbf{Z}$-points of a general linear group (see Theorem 5.2i of [22]). In certain cases (for example, in the case of potentially good reduction), the general linear group is trivial and the group embeds naturally into the $\mathbf{Q}_\ell$-points of a symplectic group. In such cases the image of the embedding lands in the $\mathbf{Z}_\ell$-points of the symplectic group, as long as $\ell$ does not divide the order of the finite inertia group or the degree of some polarization (see Theorem 5.2ii of [22]). In Theorem 5.2iii of [22] we showed that the finite inertia group is conjugate in the general linear group to a subgroup of the $\mathbf{Z}_\ell$-points of the symplectic group, even if one knows only that $\ell$ does not divide the order of the inertia group (see also Theorem 5.3 of [22]). Assuming in addition that $\ell$ is odd, then reduction modulo $\ell$ gives an embedding of the finite group in the $\mathbf{F}_\ell$-points of the symplectic group, where $\mathbf{F}_\ell$ is the finite field with $\ell$ elements. These results provide restrictions on what the finite inertia groups can be.

This leads naturally to the question of when a finite subgroup $G$ of $\mathrm{Sp}_{2d}(\mathbf{Q}_\ell)$ is conjugate in the general linear group $\mathrm{GL}_{2d}(\mathbf{Q}_\ell)$ to a subgroup of $\mathrm{Sp}_{2d}(\mathbf{Z}_\ell)$ (or can be embedded in $\mathrm{Sp}_{2d}(\mathbf{F}_\ell)$). (Of course, $G$ lies in some maximal compact subgroup of $\mathrm{Sp}_{2d}(\mathbf{Q}_\ell)$, but not necessarily in a hyperspecial one.) We treat this question in §4 and §6. Theorem 6.2 gives examples of subgroups of $\mathrm{Sp}_{2d}(\mathbf{Q}_\ell)$ of inertia type which are not conjugate in $\mathrm{GL}_{2d}(\mathbf{Q}_\ell)$ to a subgroup of $\mathrm{Sp}_{2d}(\mathbf{Z}_\ell)$. In these examples, $\ell \leq d+1$. However, it follows from Theorem 4.3 that if $\ell > d+1$, and $G$ is a subgroup of $\mathrm{Sp}_{2d}(\mathbf{Q}_\ell)$ of inertia type, then $G$ is conjugate in $\mathrm{GL}_{2d}(\mathbf{Q}_\ell)$ to a subgroup of $\mathrm{Sp}_{2d}(\mathbf{Z}_\ell)$. The results in §6 show that our bounds are sharp.

It is well-known that every finite subgroup of $\mathrm{GL}_d(\mathbf{Q}_\ell)$ is conjugate to a subgroup of $\mathrm{GL}_d(\mathbf{Z}_\ell)$ (see Part II, Chapter IV, Appendix 1 of [18]). The results in §6 show that this does not remain true in general if one replaces general linear groups by symplectic groups. However, every finite $\ell'$-subgroup of the symplectic group $\mathrm{Sp}_{2d}(\mathbf{Q}_\ell)$ is conjugate in $\mathrm{GL}_{2d}(\mathbf{Q}_\ell)$ to a subgroup of $\mathrm{Sp}_{2d}(\mathbf{Z}_\ell)$ (see Proposition 3.3 of [22]). It is known that every finite subgroup of $\mathrm{Sp}_{2d}(\mathbf{Q}_\ell)$ is conjugate in $\mathrm{Sp}_{2d}(F)$ to a subgroup of $\mathrm{Sp}_{2d}(\mathcal{O}_F)$, for some totally ramified extension $F$ of $\mathbf{Q}_\ell$ with valuation ring $\mathcal{O}_F$ (this is a special case of a general result about reductive group schemes — see Proposition 8 of [19]). We refer to [2] and [3] for a concrete description of classical groups over local fields (see also Chapters 19 and 20 in [7]).

Despite the fact that $G$ can fail to be conjugate in $\mathrm{GL}_{2d}(\mathbf{Q}_\ell)$ to a subgroup of $\mathrm{Sp}_{2d}(\mathbf{Z}_\ell)$, we prove (see Theorem 5.1) that it can nevertheless be embedded in $\mathrm{Sp}_{2d}(\mathbf{F}_\ell)$ in such a way that the characteristic polynomials are preserved (mod $\ell$), as long as $\ell > 3$. We show that these results hold for arbitrary self-dual $\ell$-adic representations, not just symplectic ones. They also hold for arbitrary finite groups, even if they are not of inertia type. Theorem 6.5 shows that the bound is sharp.

In §7 we apply the results of the earlier sections to abelian varieties, and obtain isogeny classes of abelian varieties all of whose polarizations have degree divisible by a given prime number. We give the result in the more general context of invertible sheaves.

As stated above, this paper arose from the consideration of questions in the theory of abelian varieties, and deals with questions about conjugacy of inertia subgroups of $\mathrm{Sp}_{2d}$ in $\mathrm{GL}_{2d}$ over local fields. However, this conjugacy problem may be treated as a special case of the following general problem. Given a reductive group $\mathfrak{G}$ over a local field $K$, and given a reductive subgroup $\mathfrak{H}$ of $\mathfrak{G}$, find conditions under which a finite inertia group $G \subset \mathfrak{H}(K)$ is conjugate in $\mathfrak{G}(K)$ to a subgroup of a hyperspecial subgroup of $\mathfrak{H}(K)$. (In this paper we treat the case $\mathfrak{G} = \mathrm{GL}_{2d}$, $\mathfrak{H} = \mathrm{Sp}_{2d}$.) It seems to us that this general question deserves to be studied.

Silverberg would like to thank NSA, NSF, the Science Scholars Fellowship Program at the Bunting Institute, UC Berkeley, and MSRI; Zarhin would like to thank NSF.

## 2. Definitions and notation

Suppose $\ell$ is a prime number, and $K$ is a discrete valuation field of characteristic zero and residue characteristic $\ell$. We will let $v : K^\times \to \mathbf{Z}$ denote the valuation map, normalized so that $v(K^\times) = \mathbf{Z}$. Let $\mathcal{O}_K$ denote the valuation ring, and let $e(K)$ denote the (absolute) ramification degree, i.e., $e(K) = v(\ell)$. Let $\zeta_n$ denote a primitive $n$-th root of unity in an algebraic closure of $K$, and let $\boldsymbol{\mu}_\ell$ denote the multiplicative subgroup generated by $\zeta_\ell$. If $K$ is a complete discrete valuation field, then $K$ contains $\mathbf{Q}_\ell$. If $F$ is a field and $\mathbf{Q}_\ell \subseteq F \subseteq K$, let $e(K/F)$ denote the (relative) ramification degree of $K$ over $F$.

**Definition 2.1.** A semisimple algebra is *quasisplit* if it is a direct sum of matrix algebras over (commutative) fields.

The terminology "split semisimple" is used instead of "quasisplit" in [4] (Vol. I, Defn. 3.35).

We fix once and for all a prime number $p$.

**Definition 2.2.** We say $G$ is a *group of inertia type* if $G$ is a finite group, and there exists a normal $p$-subgroup $H$ of $G$ such that the quotient $G/H$ is cyclic of order prime to $p$.

**Definition 2.3.** We say $G$ satisfies $(D_\ell)$ if $\ell$ is a prime number, $G$ is a finite group, and there is a normal subgroup $N$ of $G$ of order prime to $\ell$ such that the quotient $G/N$ is a cyclic $\ell$-group.

**Definition 2.4.** If $K$ is a field, $V$ is a $K$-vector space, $G$ is a group, and $V$ is a $G$-module, then we say $V$ is a *self-dual $K[G]$-module* if there exists a nondegenerate $G$-invariant bilinear form $f : V \times V \to K$. We say $V$ is a *symplectic* (respectively, *orthogonal*) $K[G]$-module if $f$ is alternating (respectively, symmetric).

**Definition 2.5.** Suppose $R$ is a principal ideal domain, and $T_1$ and $T_2$ are free $R$-modules. A bilinear form $f : T_1 \times T_2 \to R$ is *perfect* if the natural homomorphisms $T_1 \to \mathrm{Hom}(T_2, R)$ and $T_2 \to \mathrm{Hom}(T_1, R)$ are bijective.

**Remark 2.6.** Suppose $R$ is a principal ideal domain, $T$ is a free $R$-module of rank $2n$, and $f : T \times T \to R$ is an alternating bilinear form. If $f$ is perfect, then $\mathrm{Aut}(T, f) \cong \mathrm{Sp}_{2n}(R)$, where $\mathrm{Sp}_{2n}(R)$ denotes the group of $2n \times 2n$ symplectic matrices over $R$.

**Definition 2.7.** The *generalized quaternion group* $Q_m$ is the group of order $4m$ with the presentation:
$$Q_m = \langle a, b : a^{2m} = 1, b^2 = a^m, bab^{-1} = a^{-1} \rangle$$
(see 1.24ii, Vol. I of [4]). Note that $Q_2$ is the quaternion group of order 8.

**Definition 2.8.** Let $C_n$ denote the cyclic group of order $n$. If $p$ is odd, let $N_p$ denote the semi-direct product of $C_{2(p-1)}$ by $\boldsymbol{\mu}_p$, with the obvious action: $C_{2(p-1)} \twoheadrightarrow C_{p-1} = \mathrm{Aut}(\boldsymbol{\mu}_p)$. (Note that $N_3 = Q_3$.) Let $N_2 = Q_2$.

## 3. Lemmas

Our main technical tool in §4 is the following theorem of Serre.

**Theorem 3.1** (Serre [15]). *If $G$ is a group of inertia type, $\ell$ is a prime, and $\ell \neq p$, then $\mathbf{Q}_\ell[G]$ is quasisplit.*

**Lemma 3.2.** *If $G$ is a group of inertia type, then $G$ satisfies $(D_\ell)$ for every prime $\ell \neq p$.*

*Proof.* We use the notation of Definitions 2.2 and 2.3. Let $C = G/H$ and let $\mathcal{L}$ be the quotient of $C$ by its prime-to-$\ell$ part. Let $N$ be the kernel of the composition $G \twoheadrightarrow C \twoheadrightarrow \mathcal{L}$. Then $N$ is a normal $\ell'$-subgroup of $G$, $G/N \cong \mathcal{L}$, and $\mathcal{L}$ is a cyclic $\ell$-group. □

The next result follows immediately from the Schur-Zassenhaus Theorem.

**Lemma 3.3.** *If $G$ satisfies $(D_\ell)$, then there exists a cyclic $\ell$-subgroup $L$ of $G$ such that $G$ is the semidirect product of $N$ by $L$.*

The following result is an easy exercise.

**Lemma 3.4.** *If $G$ is a group of inertia type, and $C$ is a normal subgroup of $G$, then $G/C$ is a group of inertia type.*

**Lemma 3.5.** *Suppose $K$ is a complete discrete valuation field of characteristic zero and finite residue characteristic, $\mathfrak{m}$ is the maximal ideal of the valuation ring $\mathcal{O}_K$, and $k = \mathcal{O}_K/\mathfrak{m}$ is the residue field. Suppose $N$ is a finite group whose order is not divisible by $\mathrm{char}(k)$. Suppose $V$ is a finite-dimensional $K$-vector space and a simple $K[N]$-module, and $T$ is an $N$-stable $\mathcal{O}_K$-lattice in $V$. Then $T/\mathfrak{m}T$ is simple as a $k[N]$-module.*

*Proof.* See §14.4 and §15.5 of [17]. □

**Proposition 3.6.** *Suppose $\ell$ is a prime number, $N$ is a finite group of order prime to $\ell$, and $K$ is a complete discrete valuation field of characteristic zero and residue characteristic $\ell$. Suppose $V$ is a finite-dimensional $K$-vector space and a simple $K[N]$-module, suppose $f : V \times V \to K$ is a $K$-bilinear nondegenerate $N$-invariant pairing, suppose $T$ is an $N$-stable $\mathcal{O}_K$-lattice in $V$, and suppose $f(T,T) = \mathcal{O}_K$. Let $\mathfrak{m}$ denote the maximal ideal of $\mathcal{O}_K$. Then the restriction $f : T \times T \to \mathcal{O}_K$ is perfect, and the reduction*

$$\bar{f} : T/\mathfrak{m}T \times T/\mathfrak{m}T \to \mathcal{O}_K/\mathfrak{m}$$

*is nondegenerate.*

*Proof.* Since $V$ is a simple $K[N]$-module, $T/\mathfrak{m}T$ is a simple $(\mathcal{O}_K/\mathfrak{m})[N]$-module, by Lemma 3.5. Let $W$ denote the (left) kernel of $\bar{f}$. Since $W$ is $N$-stable, and $T/\mathfrak{m}T$ is a simple $N$-module, we have $W = 0$. Therefore, $\bar{f}$ is nondegenerate. By Nakayama's Lemma, $f : T \times T \to \mathcal{O}_K$ is perfect. □

**Lemma 3.7.** *Suppose $\ell$ is a prime number, $d$ is a positive integer, $K$ is a complete discrete valuation field of characteristic zero and residue characteristic $\ell$, $e = e(K)$, $[K(\zeta_\ell) : K] = 2d$, and $\ell = 2de + 1$. Let $T = \mathcal{O}_K[\zeta_\ell]$. Then there exists a perfect alternating $\boldsymbol{\mu}_\ell \times \{\pm 1\}$-invariant pairing*

$$f' : T \times T \to \mathcal{O}_K.$$

*Proof.* Note that $e(\mathbf{Q}_\ell(\zeta_\ell)) = \ell - 1$. Under our hypotheses, $e(K(\zeta_\ell)) = 2de = \ell - 1$, and $K(\zeta_\ell)/K$ is totally and tamely ramified. Let $\eta = \zeta_\ell - \zeta_\ell^{-1}$. Then $\eta$ is a uniformizer for $K(\zeta_\ell)$. Let $y \mapsto \bar{y}$ be the (unique) element of order two in $\mathrm{Gal}(K(\zeta_\ell)/K)$, let $\eta_1$ be a uniformizer for $K$, and define $f'$ by

$$f'(x,y) = \eta_1^{-1} \mathrm{tr}_{K(\zeta_\ell)/K}(x\eta\bar{y}).$$

Since $\ell$ does not divide $2d = [K(\zeta_\ell) : K]$, we have $\mathrm{tr}_{K(\zeta_\ell)/K}(T) = \mathcal{O}_K$. Therefore for $i = 1, \ldots, 2d$ we have

$$\eta_1 \mathcal{O}_K = \mathrm{tr}_{K(\zeta_\ell)/K}(\eta_1 T) \subseteq \mathrm{tr}_{K(\zeta_\ell)/K}(\eta^i T) \subseteq \eta^i T \cap \mathcal{O}_K \subseteq \eta_1 \mathcal{O}_K.$$

If $x \in T - \eta_1 T$, then $x = \eta^i u$ with $0 \leq i < 2d$ and $u \in T^\times$. Therefore, $f'(x,T) = \mathcal{O}_K$. It follows that $f'$ is perfect. □

**Proposition 3.8.** *Suppose $G$ is a finite group, $N$ is a normal subgroup of $G$, $W$ is a finite-dimensional vector space over a field $K$, $K[N]$ is quasisplit, and $r$ is a positive integer. We can view $\mathrm{Aut}(W)$ as naturally contained in $\mathrm{Aut}(W^r)$. Suppose $\rho : G \to \mathrm{Aut}(W^r)$ is a representation such that $\rho(N) \subset \mathrm{Aut}(W) \subseteq \mathrm{Aut}(W^r)$ and such that the restriction of $\rho$ to $N$ is irreducible. Let $Ad : G \to \mathrm{Aut}(\mathrm{End}(W^r))$ denote the corresponding adjoint representation defined by*

$$Ad(g)(u) = \rho(g)u\rho(g)^{-1}$$

*for $u \in \mathrm{End}(W^r)$ and $g \in G$. Let $E = \mathrm{End}_N(W) \subseteq \mathrm{End}(W) \subseteq \mathrm{End}(W^r)$ and let $J$ denote the image of the natural map $K[N] \to \mathrm{End}(W) \subseteq \mathrm{End}(W^r)$. Then:*

(a) $J = \mathrm{End}_E(W)$,
(b) $E = \mathrm{End}_J(W)$,
(c) $E$ is a field,
(d) the center of $J$ is $E$,
(e) $J$ and $E$ are stable under the action of $Ad(g)$ for every $g \in G$,
(f) for every $h \in N$, the action of $Ad(h)$ on $E$ is trivial.

*Proof.* By the Jacobson density theorem we have (a) and (b). Since $K[N]$ is quasisplit and $W$ is a simple $K[N]$-module, we have that $E$ is a field and $E$ is the center of $J$. Suppose $g \in G$. Since $N$ is normal in $G$, therefore $J$ is stable under $Ad(g)$. Since $E$ is the center of $J$, it follows that $E$ is stable under $Ad(g)$. Since $E$ commutes with $\rho(N)$, the action of $Ad(h)$ on $E$ is trivial for $h \in N$. □

**Corollary 3.9.** *Suppose $\ell$ is a prime number, suppose $N$ is a finite group of order prime to $\ell$, suppose $W$ is a finite-dimensional vector space over a complete discrete valuation field $K$ of characteristic zero and residue characteristic $\ell$, suppose $W$ is a simple $K[N]$-module, and suppose $K[N]$ is quasisplit. Let $E = \mathrm{End}_N(W)$. Then $E$ is a field and $E/K$ is an unramified extension.*

*Proof.* Note that if $A$ is a normal subgroup of $N$, then $K[N/A]$ is also quasisplit. Therefore, replacing $N$ by its image in $\mathrm{Aut}(W)$ we may assume that $W$ is a faithful $K[N]$-module. By Proposition 3.8, $E$ is a field. By Theorem 24.7 of [5] (see also Theorem 74.5ii, Vol. II of [4]), the field $E$ is generated over $K$ by the values of the character of the representation of $N$ on the $E$-vector space $W$. Thus, $E \subseteq K(\zeta_n)$, where $n = \#N$. Since $n$ is not divisible by $\ell$, the extension $E/K$ is unramified. □

**Lemma 3.10.** *Suppose $\ell$ is a prime number, $K$ is a complete discrete valuation field of characteristic zero and residue characteristic $\ell$, and $e = e(K)$. Suppose $A \in \mathrm{GL}_m(K)$, suppose $A$ is not a scalar, and suppose $A^\ell$ is a scalar. Then $m \geq (\ell-1)/e$.*

*Proof.* If $A^\ell$ is the scalar $c \in K^\times$, let $g(x) = x^\ell - c$. Either $g$ has a root in $K$, or $g$ is irreducible over $K$ (see the end of p. 297 of [11]). If $g$ is irreducible over $K$, then $m \geq \ell$ since $g(A) = 0$. If $g$ has a root $\gamma \in K$, then $A\gamma^{-1}$ is an element of $\mathrm{GL}_m(K)$ of exact order $\ell$. Therefore $m \geq [K(\zeta_\ell) : K] \geq (\ell-1)/e$. □

**Lemma 3.11.** *Suppose $\ell$ is an odd prime number, $K$ is an unramified extension of $\mathbf{Q}_\ell$, $M = K(\zeta_\ell)$, $\delta \in M^+ = K(\zeta_\ell + \zeta_\ell^{-1})$, $\eta = \zeta_\ell - \zeta_\ell^{-1}$, and $\mathrm{tr}_{M/K}(\delta\mathcal{O}_M) \subseteq \mathcal{O}_K$. Then*

$$\mathrm{tr}_{M/K}(\delta\eta^{\ell-2}\mathcal{O}_M) \subseteq \ell\mathcal{O}_K.$$

*Proof.* Note that $\eta$ is a uniformizer for $M$ and $\eta^2$ is a uniformizer for $M^+$. Since $e(M) = \ell - 1$, the inverse different for the extension $M/K$ is $\eta^{2-\ell}\mathcal{O}_M$ (see Proposition 13 in Chapitre III, §6 of [16]). Therefore, $\delta \in \eta^{2-\ell}\mathcal{O}_M \cap M^+ = \eta^{3-\ell}\mathcal{O}_{M^+}$. Thus,
$$\operatorname{tr}_{M/K}(\delta\eta^{\ell-2}\mathcal{O}_M) \subseteq \operatorname{tr}_{M/K}(\eta\mathcal{O}_M) = \ell\mathcal{O}_K.$$
□

The following "rigidity" result is in the spirit of similar results by Minkowski and Serre (see for example Proposition 14, III.7.6 of [1]).

**Proposition 3.12.** *Suppose $\ell$ is a prime number, and $K$ is a discrete valuation field of characteristic zero and residue characteristic $\ell$. Let $\mathfrak{m}$ denote the maximal ideal of $\mathcal{O}_K$, and let $e = e(K)$. Suppose $S$ is a free $\mathcal{O}_K$-module of finite rank, and $A$ is an automorphism of $S$ of finite order. If either*
  (a) $2e < \ell - 1$ *and* $(A - 1)^2 \in \mathfrak{m}\operatorname{End}(S)$, *or*
  (b) $e < \ell - 1$ *and* $A - 1 \in \mathfrak{m}\operatorname{End}(S)$,
*then $A = 1$.*

*Proof.* This follows directly from Theorem 6.2 of [21] for $n = \ell$, with $k = 2e$ in case (a), and with $k = e$ in case (b). □

Note that the hypothesis $2e < \ell - 1$ is satisfied if $e = 1$ and $\ell \geq 5$.

**Lemma 3.13.** *Suppose $L$ is a field of characteristic not equal to $2$, and*
$$f : L^2 \times L^2 \to L$$
*is a nondegenerate symmetric pairing. Let $g = \begin{pmatrix} 1 & 1 \\ 0 & 1 \end{pmatrix} \in \operatorname{SL}_2(L)$. Then $f$ is not $g$-invariant.*

*Proof.* Let $\{u, v\}$ denote the standard basis of $L^2$ over $L$. Suppose $f$ is $g$-invariant. Then $f(u, v) = f(gu, gv) = f(u, u+v)$, i.e., $f(u, u) = 0$. Also, $f(v, v) = f(gv, gv) = f(u+v, u+v) = f(v, v) + 2f(u, v)$, since $f$ is symmetric. Since $\operatorname{char}(L) \neq 2$, we have $f(u, v) = 0$. Therefore $f(u, w) = 0$ for every $w \in L^2$, contradicting the nondegeneracy of $f$. □

**Lemma 3.14.** *Suppose $p$ is a prime number, and $n$ is a positive integer relatively prime to $p(p-1)$. Suppose $F$ is a perfect field containing $\mathbf{F}_{p^2}$ and a primitive $n$-th root of unity. Then there exists a Galois extension $L$ of the function field $F(t)$, totally ramified at $t = 0$, such that $\operatorname{Gal}(L/F(t)) \cong N_p \times \boldsymbol{\mu}_n$. In particular, the completion of $L$ at $t = 0$ is a totally ramified Galois extension of the field $F((t))$ of formal Laurent series, with Galois group $N_p \times \boldsymbol{\mu}_n$.*

*Proof.* First suppose $p$ is odd. Let $K$ be the function field (over $F$) of the curve $y^2 = x^p - x$. The extension $K/F(x)$ is a quadratic extension which is totally ramified at $x = \infty$. One can view $N_p$ as a subgroup of $\operatorname{Aut}(K/F)$, via the automorphisms
$$(x, y) \mapsto (ax + b, cy)$$
with $a \in \mathbf{F}_p^*$, $b \in \mathbf{F}_p$, $c \in \mathbf{F}_{p^2}^*$, and $a = c^2$. Let $t = \frac{1}{(x^p - x)^{p-1}}$.

If $p = 2$, let $K$ be the function field (over $F$) of the supersingular elliptic curve in characteristic 2, $E : y^2 - y = x^3$. Note (see [23]) that $\operatorname{Aut}(E)(\cong \operatorname{SL}_2(\mathbf{F}_3)) \supset N_2 = Q_2$, with $N_2$ corresponding to the subgroup of $\operatorname{Aut}(E)$ generated by
$$(x, y) \mapsto (x + a, y + a^2 x + a)$$

for $a \in \mathbf{F}_4 - \mathbf{F}_2$. Let $t = \frac{1}{x^4-x}$.

In both cases, $N_p$ acts trivially on $F(t)$, and

$$[K : F(t)] = 2[F(x) : F(t)] = 2[F(x) : F(\frac{1}{t})] = \#N_p.$$

Therefore, $\mathrm{Gal}(K/F(t)) \cong N_p$. The extension $K/F(t)$ is totally ramified at $t = 0$. Let $L = K(t^{1/n})$. Note that $K$ and $F(t^{1/n})$ are totally ramified extensions of $F(t)$ of relatively prime degree. Therefore $L$ is a totally ramified extension of $F$, and

$$\mathrm{Gal}(L/F(t)) \cong \mathrm{Gal}(K/F(t)) \times \mathrm{Gal}(F(t^{1/n})/F(t)) \cong N_p \times \boldsymbol{\mu}_n.$$

$\square$

(See also Theorem 1 in §2 of [13], and [10].)

## 4. Embeddings over $\mathbf{Z}_\ell$

**Theorem 4.1.** *Suppose $G$ is a group that satisfies $(D_\ell)$. Suppose $K$ is a complete discrete valuation field of characteristic zero and residue characteristic $\ell$, and $e = e(K)$. Suppose $V$ is a $K$-vector space of even dimension $2d$, suppose $\ell > de+1$, suppose $V$ is a faithful simple self-dual $K[G]$-module, and suppose $K[N]$ is quasi-split. Then either:*

(a) *$V$ is simple as a $K[N]$-module, or*
(b) *$\ell = 2de + 1$, $V \cong K(\zeta_\ell)$, and either $G \cong \boldsymbol{\mu}_\ell$ or $G \cong \boldsymbol{\mu}_\ell \times \{\pm 1\}$, with the natural action on $V$.*

*Proof.* By Lemma 3.3, $G$ has a cyclic $\ell$-subgroup $L$ such that $G$ is the semidirect product of $N$ by $L$. We have

$$[K(\zeta_{\ell^2}) : K] \geq e(K(\zeta_{\ell^2})/K) \geq \ell(\ell-1)/e > \ell d > 2d.$$

Since the characteristic polynomials of the elements of $G$ acting on $V$ have degree $2d$, and $[K(\zeta_{\ell^2}) : K] > 2d$, therefore $\#L < \ell^2$. If $\#L = 1$, then $G = N$ and case (a) holds. Therefore we may assume $\#L = \ell$, and may identify $L$ with $\boldsymbol{\mu}_\ell$.

An isotypic $K[N]$-module is, by definition, a direct sum of isomorphic simple $K[N]$-modules. Let $V = \oplus_{i=1}^s V_i$ be the canonical decomposition of the $K[N]$-module $V$ into a direct sum of isotypic components (see 2.6 of [17]). The action of $G$ on $V$ induces an action of $G$ on the set of $V_i$'s. Since $V$ is a simple $K[G]$-module, the action of $G$ on $\{V_i\}$ is transitive. Since $N$ is normal in $G$, this action factors through $G/N$. Since $G/N$ has prime order $\ell$, either $s = 1$ and $V$ is isotypic, or $s = \ell$ and $V = \oplus_{g \in L}(gV_1)$. In the latter case, $2d$ is divisible by $\ell$, contradicting the hypothesis that $\ell > de + 1$. Therefore $V$ is isotypic, i.e., for some $r \in \mathbf{Z}^+$ we have $V \cong W^r$, where $W$ is a simple $K[N]$-module. The above proof followed the method of proof of Proposition 24 in §8.1 of [17]; see also Clifford's Theorem [4], Vol. 1, 11.1.

Since $V$ is a faithful $K[G]$-module, it is a faithful $K[N]$-module. Therefore, $W$ is a faithful $K[N]$-module. If $r = 1$, then (a) holds. Thus, assume $r > 1$.

Suppose $\dim_K(W) = 1$. Since $W$ is a faithful $K[N]$-module, we have $N \subseteq K^\times$. Since $V$ is a self-dual $K[N]$-module, we conclude that $N \subseteq \{\pm 1\}$. Therefore, $G$ is isomorphic to $\boldsymbol{\mu}_\ell$ or to $\boldsymbol{\mu}_\ell \times \{\pm 1\}$. Since $V$ is a simple faithful $K[G]$-module, we have $V \cong K(\zeta_\ell)$. Therefore $[K(\zeta_\ell) : K] = \dim_K(V) = 2d$, and $\ell - 1$ divides $2de$. Since $\ell - 1 > de$, therefore $\ell - 1 = 2de$ and (b) holds.

Therefore, we may assume that $\dim_K(W) \geq 2$ and $r \geq 2$. We will show that this leads to a contradiction.

Let $E = \operatorname{End}_N(W)$. We will show that $E = K$. By Corollary 3.9, $E$ is a field and $E/K$ is an unramified extension. Let $J$ denote the image of $K[N] \to \operatorname{End}(W)$. By Proposition 3.8, $J = \operatorname{End}_E(W)$, $G$ acts on $J$ and $E$ by conjugation, and the induced action of $N$ on $E$ is trivial. If $L$ acts nontrivially on $E$, then it acts faithfully, so $[E:K]$ is divisible by $\ell$. Since $\dim_K(W)$ is divisible by $[E:K]$, $r\ell$ divides $\dim_K(V) = 2d$. This contradicts the fact that $r\ell \geq 2\ell > 2(de+1) > 2d$. Therefore, $G$ commutes with $E$. Let $\gamma$ be a generator of $L$. Then $\gamma$ is an $E$-linear automorphism of $V$ of exact order $\ell$, so $\dim_E(V) \geq [E(\zeta_\ell) : E]$. Since $E/K$ is unramified,

$$\dim_E(V) \geq [E(\zeta_\ell) : E] \geq e(K(\zeta_\ell)/K) \geq (\ell-1)/e > d.$$

If $E \neq K$, then $\dim_E(V) \leq d$, a contradiction. Therefore, $E = K$.

Since $J = \operatorname{End}(W)$, therefore $G$ acts on $\operatorname{End}(W)$ by "conjugation". Since $G$ acts trivially on $E$, by the Skolem-Noether Theorem there exists $A \in \operatorname{Aut}(W) \subset \operatorname{Aut}(V)$ such that $AuA^{-1} = \gamma u \gamma^{-1}$ for every $u \in \operatorname{End}(W)$. Since $\gamma^\ell = 1$, therefore $A^\ell$ commutes with all elements of $\operatorname{End}(W)$. Thus $A^\ell$ is a scalar.

Suppose $A$ is a scalar. Then we may replace $A$ by the identity. Since $L \subseteq G \subseteq \operatorname{Aut}(V)$, we can view $\gamma$ as an element of $\operatorname{Aut}(V)$. Since $N \subset \operatorname{End}(W) \subset \operatorname{End}(V)$, we have that $\gamma$ is an element of $\operatorname{Aut}_N(V) = \operatorname{Aut}_N(W^r) \cong \operatorname{GL}_r(K)$ of exact order $\ell$. Therefore, $r \geq [K(\zeta_\ell):K] > d$. Thus, $2d = r\dim_K(W) > 2d$, a contradiction.

Therefore $A$ is not a scalar. By Lemma 3.10, we have $\dim_K(W) \geq (\ell-1)/e > d$. Therefore $2d = r\dim_K(W) > 2d$, a contradiction. $\square$

**Theorem 4.2.** *Suppose $G$ satisfies $(D_\ell)$, suppose $K$ is a complete discrete valuation field of characteristic zero and residue characteristic $\ell$, suppose $e = e(K)$, suppose $V$ is a $2d$-dimensional $K$-vector space and is a faithful simple symplectic $K[G]$-module, and suppose $\ell > de+1$. Suppose $K[N]$ is quasisplit. Then there exist a $G$-stable $\mathcal{O}_K$-lattice $T$ in $V$ and a perfect alternating $G$-invariant bilinear form*

$$f' : T \times T \to \mathcal{O}_K.$$

*Proof.* By Theorem 4.1, we are in case (a) or (b). In case (b), apply Lemma 3.7. In case (a), let $f$ be the associated alternating pairing on $V$, and let $T$ be a $G$-stable $\mathcal{O}_K$-lattice in $V$. By multiplying $f$ by a suitable power of a uniformizer, we can assume that $f(T,T) = \mathcal{O}_K$. Now apply Proposition 3.6. $\square$

**Theorem 4.3.** *Suppose $G$ is a group of inertia type, suppose $\ell$ is a prime number, and suppose $\ell \neq p$. Suppose $K$ is a complete discrete valuation field of characteristic zero and residue characteristic $\ell$, suppose $e = e(K)$, suppose $V$ is a $2d$-dimensional $K$-vector space and is a symplectic $K[G]$-module, and suppose $\ell > de+1$. Then there exist a $G$-stable $\mathcal{O}_K$-lattice $T$ in $V$ and a perfect alternating $G$-invariant bilinear form*

$$f' : T \times T \to \mathcal{O}_K.$$

*Proof.* We can reduce to the case where $V$ is a simple symplectic $K[G]$-module as follows. Let $f$ be the associated alternating pairing on $V$. Write $V$ as a direct sum of simple $K[G]$-modules

$$V = V_1 \oplus \cdots \oplus V_k.$$

If for some $i$ the restriction $f_i$ of $f$ to $V_i$ is not identically zero, then $f_i$ is nondegenerate. Let $V_i'$ be the orthogonal complement of $V_i$ in $V$ with respect to $f$. Clearly, the restriction of $f$ to $V_i'$ is also nondegenerate. Apply induction to the symplectic $K[G]$-modules $V_i$ and $V_i'$, to obtain the desired result. We may therefore assume that $f(V_i, V_i) = 0$ for every $i$. Since $f$ is nondegenerate, there exists a $V_j$ such that $f(V_1, V_j) \neq 0$. Since $V_1$ and $V_j$ are simple $K[G]$-modules, $f : V_1 \times V_j \to K$ is nondegenerate. Thus the restriction $\tilde{f}$ of $f$ to $W = V_1 \oplus V_j$ is nondegenerate. Let $V''$ denote the orthogonal complement of $W$ in $V$ with respect to $f$. Again by induction, we obtain a $G$-stable lattice $T''$ in $V''$ and a $G$-invariant perfect $\mathcal{O}_K$-valued alternating pairing $f''$ on $T''$. Let $T_1$ be a $G$-stable $\mathcal{O}_K$-lattice in $V_1$, let

$$T_j = \{x \in V_j : f(T_1, x) \subseteq \mathcal{O}_K\},$$

let $T = T_1 \oplus T_j \oplus T''$, and let $f' = \tilde{f} \oplus f''$.

By Lemma 3.4, replacing $G$ by its image in $\mathrm{Aut}(V)$, we may reduce to the case where $V$ is a faithful simple symplectic $K[G]$-module. Suppose $\ell$ is a prime and $\ell \neq p$. By Lemma 3.2, $G$ satisfies $(D_\ell)$, and we have a corresponding normal subgroup $N$. By Theorem 3.1, $\mathbf{Q}_\ell[N]$ is quasisplit. Therefore, $K[N]$ is quasisplit. The theorem now follows from Theorem 4.2. $\square$

**Remark 4.4.** If $d = 1$, then $\dim_K(V) = 2$ and $\mathrm{Aut}(V, f) = \mathrm{SL}(V)$. In this case the conclusion of Theorems 4.2 and 4.3 holds without the requirement that $\ell > de + 1$. Indeed, one can let $T$ be any $G$-stable $\mathcal{O}_K$-lattice in $V$, and let $f' : T \times T \to \mathcal{O}_K$ be any alternating perfect pairing.

**Remark 4.5.** It has been pointed out to us that our arguments allow one to drop the assumption that $K[N]$ is quasisplit from Theorem 4.1 (and therefore also from Theorem 4.2), and thereby replace "$G$ is a group of inertia type" by "$G$ satisfies $(D_\ell)$" in the statement of Theorem 4.3. This follows since one may pass to an unramified extension $K'/K$ such that $K'[N]$ is quasisplit.

## 5. Embeddings over $\mathbf{F}_\ell$

The following result holds for arbitrary finite groups (not just those of inertia type) and for arbitrary nondegenerate self-dual pairings (not just alternating ones).

**Theorem 5.1.** *Suppose $\ell$ is a prime number, and $K$ is a discrete valuation field of characteristic zero and residue characteristic $\ell$. Let $\mathfrak{m}$ denote the maximal ideal, let $k = \mathcal{O}_K/\mathfrak{m}$, let $e = e(K)$, and suppose $2e < \ell - 1$. Suppose $V$ is a $K$-vector space of finite dimension $n$, suppose $f : V \times V \to K$ is a nondegenerate alternating (respectively, symmetric) $K$-bilinear form, suppose $G$ is a finite group, and suppose*

$$\rho : G \hookrightarrow \mathrm{Aut}_K(V, f)$$

*is a faithful representation of $G$ on $V$ that preserves the form $f$. Then there exist a nondegenerate alternating (respectively, symmetric) $k$-valued $k$-bilinear form $f_0$ on $k^n$, and a faithful representation*

$$\bar\rho : G \hookrightarrow \mathrm{Aut}_k(k^n, f_0),$$

*such that for every $g \in G$, the characteristic polynomial of $\bar\rho(g)$ is the reduction modulo $\mathfrak{m}$ of the characteristic polynomial of $\rho(g)$.*

*Proof.* If $S$ is a $G$-stable $\mathcal{O}_K$-lattice in $V$, let

$$S^* = \{x \in V : f(x, S) \subseteq \mathcal{O}_K\}.$$

Fix a $G$-stable $\mathcal{O}_K$-lattice $S$ in $V$. Let $\eta$ denote a uniformizer for $\mathcal{O}_K$. Multiplying $f$ by an integral power of $\eta$ if necessary, we may assume that $f(S,S) = \mathcal{O}_K$. Then $S \subseteq S^*$. Let $S_0 = S$ and let

$$S_{i+1} = S_i + (\eta^{-1} S_i \cap \eta S_i^*) \quad \text{for} \quad i \geq 0.$$

Then $S_i$ is a $G$-stable $\mathcal{O}_K$-lattice in $V$, $f(S_i, S_i) = \mathcal{O}_K$, and $S_i \subseteq S_{i+1} \subseteq S_{i+1}^* \subseteq S_i^*$. Note that $S_{i+1} = S_i$ if and only if $\eta S_i^* \subseteq S_i$. We have

$$S = S_0 \subseteq S_1 \subseteq S_2 \subseteq \ldots \subseteq S^*.$$

Since $S^*/S$ is finite, we have $S_j = S_{j+1}$ for some $j$. Let $T = S_j$. Then $T$ is a $G$-stable $\mathcal{O}_K$-lattice in $V$ such that $f(T,T) = \mathcal{O}_K$ and $\mathfrak{m} T^* = \eta T^* \subseteq T$.

Let $\bar{f} : T/\mathfrak{m} T \times T/\mathfrak{m} T \to k$ be the reduction of $f$ modulo $\mathfrak{m}$. Then $\ker(\bar{f}) = \mathfrak{m} T^*/\mathfrak{m} T \cong T^*/T$. Clearly, $\bar{f}$ is nondegenerate on $(T/\mathfrak{m} T)/\ker(\bar{f}) \cong T/\mathfrak{m} T^*$. On $T^* \times T^*$, the form $\eta f$ is $\mathcal{O}_K$-valued. Let $\tilde{f}$ denote the reduction modulo $\mathfrak{m}$ of the restriction of $\eta f$ to $T^* \times T^*$. Since $(T^*)^* = T$, we have $\ker(\tilde{f}) = T/\mathfrak{m} T^*$. Therefore, $\tilde{f}$ is nondegenerate on $(T^*/\mathfrak{m} T^*)/(T/\mathfrak{m} T^*) \cong T^*/T \cong \ker(\bar{f})$. We thus obtain a homomorphism

$$\psi : G \to \operatorname{Aut}_k((T/\mathfrak{m} T)/\ker(\bar{f}), \bar{f}) \times \operatorname{Aut}_k(\ker(\bar{f}), \tilde{f}) \cong$$

$$\operatorname{Aut}_k(T/\mathfrak{m} T^*, \bar{f}) \times \operatorname{Aut}_k(T^*/T, \tilde{f}) \hookrightarrow \operatorname{Aut}_k(T^*/\mathfrak{m} T^*, \bar{f} \times \tilde{f}) \cong \operatorname{Aut}_k(k^n, f_0)$$

for an appropriate pairing $f_0$. All the elements $\sigma \in \ker(\psi)$ act as the identity on $T/\mathfrak{m} T^*$ and on $T^*/T$, and thus $(\sigma - 1)^2 T^* \subseteq \mathfrak{m} T^*$. Proposition 3.12a implies $\psi$ is injective. Let $\bar{\rho} = \psi$. □

**Remark 5.2.** If one is concerned only with preserving the characteristic polynomials, and does not insist that $\bar{\rho}$ be an embedding, then in the symplectic case the above result can instead be achieved with the aid of Proposition 8 of [19], rather than Proposition 3.12a above.

**Remark 5.3.** In the symplectic case, one can obtain the result by using the description of maximal bounded subgroups (maximal compact subgroups in the case of locally compact $K$) from the theory of buildings. The maximal bounded subgroups $\mathcal{G}$ of $\operatorname{Sp}_{2d}(K)$ are the stabilizers of lattices $\Lambda$ which are orthogonal direct sums

$$(\Lambda_1, f_1) \oplus (\Lambda_2, \eta f_2)$$

where $\eta$ is a uniformizer, and $f_1$ and $f_2$ are perfect alternating forms on the lattices $\Lambda_1$ and $\Lambda_2$, respectively. We may assume $G \subset \mathcal{G}$. The action of $\mathcal{G}$ on $\Lambda/\eta\Lambda$ induces the desired map

$$\bar{\rho} : G \hookrightarrow \operatorname{Sp}(\Lambda_1/\eta\Lambda_1) \times \operatorname{Sp}(\Lambda_2/\eta\Lambda_2)$$

(if $g \in \ker(\bar{\rho})$, then $(g-1)^2 \equiv 0 \pmod{\eta}$; if $g$ has finite order, then $g = 1$ under our hypotheses).

**Theorem 5.4.** *Suppose $\ell$ is a prime number, and $L$ is a discrete valuation field of characteristic zero and residue characteristic $\ell$. Suppose $K$ is a quadratic extension of $L$, let $\mathfrak{m}$ denote the maximal ideal, and let $k = \mathcal{O}_K/\mathfrak{m}$. Let $e = e(K)$, and suppose $2e < \ell - 1$. Suppose $V$ is a $K$-vector space of finite dimension $n$, suppose $f : V \times V \to K$ is a nondegenerate pairing which is hermitian (respectively, skew-hermitian) with respect to the extension $K/L$, suppose $G$ is a finite group, and suppose $\rho : G \hookrightarrow \operatorname{Aut}_K(V, f)$ is a faithful representation. Suppose $K/L$ is unramified and thus $k$ is a quadratic extension of the residue field $k_L$ of $L$. Then there exist a nondegenerate* $k$-valued pairing $f_0$ on $k^n$ which is hermitian (respectively, skew-hermitian) with respect to the extension $k/k_L$, and a faithful representation $\bar{\rho} : G \hookrightarrow \mathrm{Aut}_k(k^n, f_0)$, such that for every $g \in G$, the characteristic polynomial of $\bar{\rho}(g)$ is the reduction modulo $\mathfrak{m}$ of the characteristic polynomial of $\rho(g)$.

*Proof.* Let $\eta$ denote a uniformizer for $L$. Then $\eta$ is also a uniformizer for $K$, and the proof is a repetition of the proof of Theorem 5.1. The reductions $\bar{f}$ and $\tilde{f}$ are now hermitian (respectively, skew-hermitian). □

**Remark 5.5.** In the setting of Theorem 5.4, suppose now that $K/L$ is ramified. Then $k_L = k$, and one can write $K = L(\sqrt{D})$ where $D$ is a uniformizer for $L$ and $\eta = \sqrt{D}$ is a uniformizer for $K$. Let $S$ be a $G$-stable $\mathcal{O}_K$-lattice in $V$. Then $\eta^r f(S, S) = \mathcal{O}_K$ for some integer $r$. If $r$ is even then $\eta^r f$ is hermitian (respectively, skew-hermitian) and its reduction is symmetric (respectively, alternating). Further, $\eta^{r+1} f$ is skew-hermitian (respectively, hermitian) and its reduction is alternating (respectively, symmetric). If $r$ is odd, then $\eta^r f$ is skew-hermitian (respectively, hermitian) and its reduction is alternating (respectively, symmetric). Further, $\eta^{r+1} f$ is hermitian (respectively, skew-hermitian) and its reduction is symmetric (respectively, alternating). In all cases one can proceed as in Theorem 5.1 and obtain an embedding of $G$ into a product of an orthogonal group $\mathrm{O}_s(k)$ and a symplectic group $\mathrm{Sp}_{n-s}(k)$, which "respects" the characteristic polynomials.

## 6. Sharpness

The next proposition shows that one cannot replace "symplectic" by "orthogonal" in the statements of Theorems 4.2 and 4.3. Here, $e = 1$ and $\ell - 1 = 2d = 2de > de$.

**Proposition 6.1.** *Suppose $\ell$ is an odd prime. Let $G = \boldsymbol{\mu}_\ell$, let $M = \mathbf{Q}_\ell(\zeta_\ell)$, and let $V$ be $M$ viewed as a $\mathbf{Q}_\ell$-vector space. Then $V$ carries a natural structure of a faithful simple orthogonal $\mathbf{Q}_\ell[G]$-module, but there is no $G$-stable lattice $T$ in $V$ with a perfect symmetric $G$-invariant pairing $f' : T \times T \to \mathbf{Z}_\ell$.*

*Proof.* Define
$$f : V \times V \to \mathbf{Q}_\ell \quad \text{by} \quad f(x, y) = \mathrm{tr}_{M/\mathbf{Q}_\ell}(x\bar{y}),$$
where $y \mapsto \bar{y}$ is the automorphism of $M$ over $\mathbf{Q}_\ell$ which sends $\zeta_\ell$ to its inverse. Then $f$ is a nondegenerate symmetric $G$-invariant pairing. Let $S = \mathbf{Z}_\ell[\zeta_\ell]$ and let $M^+ = \mathbf{Q}_\ell(\zeta_\ell + \zeta_\ell^{-1})$. We have $\mathrm{End}_G(V) = M$. Suppose $T$ is a $G$-stable $\mathbf{Z}_\ell$-lattice in $V$ and $f' : T \times T \to \mathbf{Z}_\ell$ is a $G$-invariant perfect symmetric bilinear form. Every nondegenerate $G$-invariant symmetric bilinear form on $V$ is of the form $f_\delta(x, y) = f(\delta x, y)$ for some $\delta \in (M^+)^\times$. Therefore $f' = f_\delta$ for some $\delta \in M^+$. Let $\eta = \zeta_\ell - \zeta_\ell^{-1}$. Then $\eta$ is a uniformizer for $M$ and $\eta^2$ is a uniformizer for $M^+$. Since $T$ is a $G$-stable $\mathbf{Z}_\ell$-lattice in $V$, therefore $T$ is a $\mathbf{Z}_\ell[\zeta_\ell]$-lattice. Since $V$ is a one-dimensional $\mathbf{Q}_\ell(\zeta_\ell)$-vector space, it follows that $T = \eta^r S$ for some $r \in \mathbf{Z}$. Replacing $\delta$ by $\delta(\eta\bar{\eta})^r$, we can reduce to the case where $f_\delta$ is a perfect $\mathbf{Z}_\ell$-valued pairing on $S$. Since $f_\delta(S, S) \subseteq \mathbf{Z}_\ell$, therefore $\mathrm{tr}_{M/\mathbf{Q}_\ell}(\delta S) \subseteq \mathbf{Z}_\ell$. By Lemma 3.11, $f_\delta(\eta^{\ell-2} S, S) \subseteq \ell \mathbf{Z}_\ell$. The perfectness of $f_\delta : S \times S \to \mathbf{Z}_\ell$ is now contradicted by the fact that $\eta^{\ell-2} S$ contains elements of $S - \ell S$. □

The following results show that the bound $\ell > de + 1$ in §4 is sharp.

**Theorem 6.2.** *Suppose $\ell$ is an odd prime, $N$ is a finite group of order not divisible by $\ell$, $K$ is a complete discrete valuation field of characteristic zero and residue characteristic $\ell$ which is an unramified extension of $\mathbf{Q}_\ell$, and $W$ is a non-zero $K$-vector space which is also an absolutely simple faithful symplectic $K[N]$-module. Let $G = N \times \boldsymbol{\mu}_\ell$, and let $V = W \otimes_K K(\zeta_\ell)$. Then $V$ carries a natural structure of a faithful symplectic $K[G]$-module, but there is no $G$-stable $\mathcal{O}_K$-lattice $T$ in $V$ with a perfect alternating $G$-invariant pairing $f' : T \times T \to \mathcal{O}_K$.*

*Proof.* Let $T_1$ be an $N$-stable $\mathbf{Z}_\ell$-lattice in $W$. Multiplying the given alternating $N$-invariant pairing $f_1$ on $W$ by a suitable power of $\ell$, we may assume that $f_1(T_1, T_1) = \mathcal{O}_K$. Since $\#N$ is not divisible by $\ell$, it follows from Proposition 3.6 that $f_1 : T_1 \times T_1 \to \mathcal{O}_K$ is perfect.

Let $M = K(\zeta_\ell)$, and let $W_2$ be $M$ viewed as a $K$-vector space. Define

$$f_2 : W_2 \times W_2 \to K \quad \text{by} \quad f_2(\alpha, \beta) = \text{tr}_{M/K}(\alpha \bar{\beta}),$$

where $\beta \mapsto \bar{\beta}$ is the automorphism of $M$ over $K$ which sends $\zeta_\ell$ to its inverse. Then $f_2$ is a nondegenerate $\boldsymbol{\mu}_\ell$-invariant symmetric pairing. Let $f = f_1 \otimes f_2$. Then $f$ is a nondegenerate $G$-invariant alternating pairing, with respect to which $V$ is a faithful symplectic $K[G]$-module. We have $\text{End}_N(W) = K$. Since $\text{End}_{\boldsymbol{\mu}_\ell}(W_2) = M$, we have $\text{End}_G(V) = M$. It therefore follows from the definition of $f$ that every nondegenerate $G$-invariant alternating bilinear form on $V$ is of the form $f_\delta(x, y) = f(\delta x, y)$ for some $\delta \in (M^+)^\times$, where $M^+ = K(\zeta_\ell + \zeta_\ell^{-1})$. Let $S = T_1 \otimes_{\mathcal{O}_K} \mathcal{O}_M \subset V$, and let $k$ denote the residue field of $K$. By Lemma 3.5, $S/(1 - \zeta_\ell)S$ ($= T_1/\ell T_1$) is a simple $k[N]$-module, and therefore is a simple $k[G]$-module. Let $\eta = \zeta_\ell - \zeta_\ell^{-1}$. Every $G$-stable $\mathcal{O}_K$-lattice in the $K$-vector space $V$ is an $\mathcal{O}_M$-lattice. By Exercise 15.3 of [17], it follows that every $G$-stable lattice in $V$ is of the form $\eta^r S$ with $r \in \mathbf{Z}$. Thus if $T$ is a $G$-stable lattice in $V$ and $f' : T \times T \to \mathbf{Z}_\ell$ is a $G$-invariant perfect alternating pairing, then there exist $r \in \mathbf{Z}$ and $\delta \in M^+$ such that $T = \eta^r S$ and $f' = f_\delta$. Replacing $\delta$ by $\delta(\eta \bar{\eta})^r$, we can reduce to the case where $f_\delta$ is a perfect $\mathcal{O}_K$-valued pairing on $S$. Since $f_\delta(S, S) \subseteq \mathcal{O}_K$, therefore $\text{tr}_{M/K}(\delta S) \subseteq \mathcal{O}_K$. By Lemma 3.11, $f_\delta(\eta^{\ell-2} S, S) \subseteq \ell \mathcal{O}_K$. The perfectness of $f_\delta : S \times S \to \mathcal{O}_K$ is contradicted by the fact that $\eta^{\ell-2} S$ contains elements of $S - \ell S$. $\square$

Applying Theorem 6.2 to the family of examples in the next result, we obtain examples which show that for all primes $p$ and for all odd primes $\ell \neq p$, the bound $\ell > de + 1$ in Theorem 4.3 is sharp. In these examples we have that $e = 1$ and $\ell - 1 = d = de$, and the group $G = N \times \boldsymbol{\mu}_\ell$ is a group of inertia type.

**Proposition 6.3.** *Let $q = p$ if $p$ is odd, and let $q = 4$ if $p = 2$. Assume that the quadruple $(\ell, K, W, N)$ is as in one of the following two cases.*
  (a) *Let $\ell$ be an odd prime different from $p$, let $K = \mathbf{Q}_\ell(\zeta_q)$, let $W = K^2$, and let $N$ be the subgroup of $\text{SL}(W)$ generated by $\begin{pmatrix} \zeta_q & 0 \\ 0 & \zeta_q^{-1} \end{pmatrix}$ and $\begin{pmatrix} 0 & -1 \\ 1 & 0 \end{pmatrix}$.*
  (b) *Let $\ell$ be an odd prime different from $p$ such that $\ell \equiv -1 \pmod{q}$. Let $K = \mathbf{Q}_\ell$, and let $W$ be the 2-dimensional $\mathbf{Q}_\ell$-vector space $\mathbf{Q}_\ell(\zeta_q)$. Let $\sigma$ be the nontrivial automorphism of the quadratic extension $\mathbf{Q}_\ell(\zeta_q)/\mathbf{Q}_\ell$. Pick $\alpha \in \mathbf{Q}_\ell(\zeta_q)$ such that $\text{Norm}_{\mathbf{Q}_\ell(\zeta_q)/\mathbf{Q}_\ell}(\alpha) = -1$ (this can be done since $\mathbf{Q}_\ell(\zeta_q)/\mathbf{Q}_\ell$ is unramified). Let $N$ be the subgroup of $\text{SL}(W)$ generated by $\zeta_q$ (with the natural action) and the element $\tau$ of order 4 defined by $\tau(x) = \alpha \sigma(x)$ for $x \in W$.*

Then $\ell$, $K$, $N$, and $W$ satisfy the hypotheses of Theorem 6.2, and $G = N \times \boldsymbol{\mu}_\ell$ is a group of inertia type.

*Proof.* The group $N$ is the generalized quaternion group $Q_p$ of order $4p$. Since $N$ is non-abelian and $W$ is a 2-dimensional faithful $K[N]$-module, it follows that $W$ is an absolutely simple $K[N]$-module. If $p$ is odd, then $\boldsymbol{\mu}_p$ is a normal subgroup of $N$ and $G$, and $G/\boldsymbol{\mu}_p$ is cyclic of order $4\ell$. If $p = 2$, then $N$ is a normal Sylow-2-subgroup of $G$, and $G/N$ is cyclic of order $\ell$. In each case, $G$ is of inertia type. $\square$

Note that $\mathbf{Q}_\ell(\zeta_q) = \mathbf{Q}_\ell$ if and only if $\ell \equiv 1 \pmod{q}$. Therefore when $p = 2$ or $3$, the preceding result yields examples which establish the sharpness of the bound $\ell > de + 1$ in Theorem 4.3 with $K = \mathbf{Q}_\ell$ and $d = \ell - 1$, for all odd primes $\ell \neq p$.

When $\ell$, $K$, $N$, and $W$ are as in Proposition 6.3, and $b$ is a positive integer, then the following result gives examples where $e = 1$ and $de = d = \ell - 1 + b$, thereby providing examples which do not satisfy the conclusion of Theorem 4.3, for all odd primes $\ell \leq de + 1$ different from $p$.

**Corollary 6.4.** *Suppose $\ell$, $N$, $K$, $W$, $V$, and $G$ are as in Theorem 6.2, and $b \in \mathbf{Z}^+$. Let $W_0$ be the $K$-vector space $K^{2b}$ with trivial $G$-action, and let $U$ be the $K[G]$-module $V \oplus W_0$. Then $U$ carries a natural structure of a faithful symplectic $K[G]$-module, but there is no $G$-stable $\mathcal{O}_K$-lattice $T$ in $U$ with a perfect alternating $G$-invariant pairing $f' : T \times T \to \mathcal{O}_K$.*

*Proof.* The natural symplectic structure on $U$ arises from the sum of the natural alternating pairing on $V$ and the standard alternating pairing on $W_0$. Let $\tau = \frac{1}{\#N}\sum_{h \in N} h \in \mathbf{Z}_\ell[N]$. Then $\tau^2 = \tau$, $\tau U = W_0$, and $(1-\tau)U = V$. Suppose $T$ is a $G$-stable $\mathcal{O}_K$-lattice in $U$ and $f' : T \times T \to \mathcal{O}_K$ is a perfect alternating $G$-invariant pairing. Let $T_1 = (1-\tau)T$ and let $T_2 = \tau T$. Then $T = T_1 \oplus T_2$, and $T_1$ (respectively, $T_2$) is a $G$-stable $\mathcal{O}_K$-lattice in $V$ (respectively, $W_0$). Since $f'$ is $G$-invariant and the $G$-action on $W_0$ is trivial, we have $f'(x,y) = f'(gx, gy) = f'(gx, y)$ for every $x \in T_1$, $y \in T_2$, and $g \in G$. Therefore, $f'((1-g)T_1, T_2) = 0$ for every $g \in G$. Since $V$ is a simple $K[G]$-module, therefore the $(1-g)T_1$'s generate $T_1$. Thus, $f'(T_1, T_2) = 0$. Since $f'$ is perfect, it follows that the restriction of $f'$ to $T_1 \times T_1$ is perfect, contradicting Theorem 6.2. $\square$

Theorems 6.5 and 6.6 below enable us to obtain examples with $e = (\ell-1)/2$. In particular, Theorem 6.5 gives examples which show that the bound $\ell > 2e + 1$ in Theorem 5.1 is sharp. The proof of Theorem 6.6 uses Theorem 6.5.

**Theorem 6.5.** *Suppose $\ell$ is an odd prime number, $F$ is a complete discrete valuation field of characteristic zero and residue characteristic $\ell$ which is an unramified extension of $\mathbf{Q}_\ell$, and $N$ is a finite group of order not divisible by $\ell$. Let $K = F(\zeta_\ell + \zeta_\ell^{-1})$, let $\mathfrak{m}$ denote the maximal ideal of $\mathcal{O}_K$, and let $k = \mathcal{O}_K/\mathfrak{m}$. Suppose $W$ is a non-zero $K$-vector space which is also an absolutely simple faithful symplectic $K[N]$-module, and write $\dim_K(W) = 2t$. Let $V$ be the $K$-vector space $W \otimes_K K(\zeta_\ell)$, and let $G = N \times \boldsymbol{\mu}_\ell$. Then there are a natural nondegenerate alternating $K$-bilinear form $f : V \times V \to K$ and a natural faithful irreducible representation $\rho : G \hookrightarrow \mathrm{Aut}_K(V, f)$. If $L$ is a field which contains $k$, and $f_0 : L^{4t} \times L^{4t} \to L$ is a nondegenerate $L$-bilinear alternating form, then there does not exist a faithful representation $\bar{\rho} : G \hookrightarrow \mathrm{Aut}_L(L^{4t}, f_0)$ having the property that for every $g \in G$, $\mathrm{tr}(\bar{\rho}(g)) \equiv \mathrm{tr}(\rho(g)) \pmod{\mathfrak{m}}$.*

*Proof.* Let $M = F(\zeta_\ell)$, let $x \mapsto \bar{x}$ denote the nontrivial automorphism of $M$ over $K$, and let $V_2$ be $M$ viewed as a 2-dimensional $K$-vector space. By our assumptions, we have an absolutely irreducible, faithful representation $\rho_1 : N \hookrightarrow \mathrm{Aut}_K(W, f_1)$ where $f_1$ is a nondegenerate alternating $K$-valued pairing on $W$. The natural action of $\boldsymbol{\mu}_\ell$ on $V_2$ respects the nondegenerate symmetric bilinear form $f_2 : V_2 \times V_2 \to K$ defined by $f_2(x, y) = \mathrm{tr}_{M/K}(x\bar{y})$. We thus obtain a faithful representation $\rho_2 : \boldsymbol{\mu}_\ell \hookrightarrow \mathrm{Aut}_K(V_2, f_2)$. Let $\rho$ be the representation

$$\rho : G = N \times \boldsymbol{\mu}_\ell \hookrightarrow \mathrm{Aut}_K(W, f_1) \times \mathrm{Aut}_K(V_2, f_2) \hookrightarrow \mathrm{Aut}_K(V, f)$$

defined by $\rho(a, b) = \rho_1(a) \otimes \rho_2(b)$, where $f = f_1 \otimes f_2$.

Suppose $L$ is a field which contains $k$, $V_0$ is a $4t$-dimensional $L$-vector space, $f_0$ is a nondegenerate $L$-valued alternating pairing on $V_0$, and $\bar{\rho} : G \hookrightarrow \mathrm{Aut}_L(V_0, f_0)$ is a faithful representation. Since $\#N$ is not divisible by $\ell$, $V_0$ is a semisimple $L[N]$-module. Let $T$ be an $N$-stable $\mathcal{O}_K$-lattice $T$ in $W$. By Proposition 3.6, for some non-zero multiple $f'$ of $f_1$ we have $\rho_1(N) \subset \mathrm{Aut}(T, f')$, and the reduction $\bar{f}'$ is nondegenerate. The composition of $\rho_1$ with the reduction map gives a faithful representation $\bar{\rho}_1 : N \hookrightarrow \mathrm{Aut}(W_0, \bar{f}')$, where the corresponding $k[N]$-module $W_0$ is absolutely simple and symplectic. By Schur's Lemma, every $N$-invariant bilinear form on $W_0$ is alternating. If $\mathrm{tr}((\bar{\rho}_1 \oplus \bar{\rho}_1)(h)) = \mathrm{tr}(\bar{\rho}(h))$ for every $h \in N$, then the multiplicities of the simple components of the $L[N]$-modules $V_0$ and $W_0 \oplus W_0$ are congruent modulo $\ell$ (see Theorem 17.3 of [4]). Then $V_0$ and $W_0^2$ are isomorphic $L[N]$-modules, since $W_0$ is simple and $\ell > 2$. Since $\mathrm{End}_N(W_0) = F$, we have $\mathrm{End}_N(V_0) = \mathrm{M}_2(F)$. Fix a generator $c$ of $\boldsymbol{\mu}_\ell$. Then $c$ is an element of $\mathrm{End}_N(V_0)$ of multiplicative order $\ell$. We can therefore identify $V_0$ with $W_0 \oplus W_0$ in such a way that $c(x, y) = (x + y, y)$ for every $(x, y) \in W_0 \oplus W_0 = V_0$. For $x, y \in W_0$ we have

$$f_0((x, 0), (0, y)) = f_0(c(x, 0), c(0, y)) = f_0((x, 0), (y, y)) =$$

$$f_0((x, 0), (0, y)) + f_0((x, 0), (y, 0)).$$

Therefore,

(1) $\qquad f_0((x, 0), (y, 0)) = 0$ for all $x, y \in W_0$.

Further,

$$f_0((0, x), (0, y)) = f_0(c(0, x), c(0, y)) = f_0((x, x), (y, y)) =$$

$$f_0((0, x), (0, y)) + f_0((x, 0), (y, 0)) + f_0((x, 0), (0, y)) + f_0((0, x), (y, 0)).$$

Therefore, $f_0((x, 0), (0, y)) + f_0((0, x), (y, 0)) = 0$. Since $f_0$ is alternating,

(2) $\qquad f_0((x, 0), (0, y)) = -f_0((0, x), (y, 0)) = f_0((y, 0), (0, x)).$

Define $h : W_0 \times W_0 \to L$ by $h(x, y) = f_0((x, 0), (0, y))$. By (2), $h$ is symmetric. Since $f_0$ is $N$-invariant, so is $h$. The nondegeneracy of $h$ follows from (1) and the nondegeneracy of $f_0$. Since $h$ is a $N$-invariant pairing on $W_0$, it is alternating. Since $h$ is both alternating and symmetric we have $h = 0$, giving a contradiction. □

The next result shows that the bound $\ell > de + 1$ in Theorem 4.3 is sharp with respect to $e$ (subject to the restrictions imposed by Remark 4.4). Here we have $e = (\ell - 1)/2$, $d = 2$, and $\ell - 1 = de$. A special case of Theorem 6.5 and 6.6 is (in the notation of Proposition 6.3) when $\ell$ is odd, $F = \mathbf{Q}_\ell(\zeta_q)$ (resp., $\mathbf{Q}_\ell$), $W = \mathbf{Q}_\ell(\zeta_q)^2$ (resp., $\mathbf{Q}_\ell(\zeta_q)$), and $N$ is the generalized quaternion group $Q_p$. The group $G$ is then of inertia type.

**Theorem 6.6.** *Suppose $\ell$ is an odd prime number, $F$ is a complete discrete valuation field of characteristic zero and residue characteristic $\ell$ which is an unramified extension of $\mathbf{Q}_\ell$, and $N$ is a finite group of order not divisible by $\ell$. Let $K = F(\zeta_\ell + \zeta_\ell^{-1})$, and suppose $W$ is a non-zero $K$-vector space which is also an absolutely simple faithful symplectic $K[N]$-module. Let $V$ be the $K$-vector space $W \otimes_K K(\zeta_\ell)$, and let $G = N \times \boldsymbol{\mu}_\ell$. Then $V$ carries a natural structure of a faithful symplectic $K[G]$-module, and there is no $G$-stable lattice $T$ in $V$ with a $G$-invariant perfect alternating pairing $f' : T \times T \to \mathcal{O}_K$.*

*Proof.* Let $\mathfrak{m}$ denote the maximal ideal of $\mathcal{O}_K$, let $k = \mathcal{O}_K/\mathfrak{m}$, let $e = e(K)$, and write $\dim_K(W) = 2t$. By Theorem 6.5, $V$ carries a natural structure of a faithful symplectic $K[G]$-module. Suppose there existed a $G$-stable lattice $T$ in $V$ and a $G$-invariant perfect alternating pairing $f' : T \times T \to \mathcal{O}_K$, and consider the resulting representation $\rho : G \hookrightarrow \mathrm{Aut}_K(T, f')$. The composition of $\rho$ with reduction modulo $\mathfrak{m}$ gives a representation $\bar\rho : G \to \mathrm{Aut}_k(k^{4t}, \bar{f}')$, where $\bar{f}'$ is the reduction of $f'$. If $A \in \ker(\bar\rho)$, then $A - 1 \in \mathfrak{m}\mathrm{End}(T)$. By Proposition 3.12b we have $A = 1$, since $e = (\ell-1)/2 < \ell - 1$. Therefore $\bar\rho$ is injective, contradicting Theorem 6.5. $\square$

## 7. Polarizations and isogeny classes of abelian varieties

In Theorem 7.1 we use the results of §6 to construct isogeny classes of abelian varieties which admit no principal polarization. In particular, for every odd prime $\ell$, and every integer $d \geq \ell - 1$, we construct isogeny classes of $d$-dimensional supersingular abelian varieties $B$ such that the degree of every polarization on $B$ is divisible by $\ell$.

Suppose $F$ is a field, let $F^s$ denote a separable closure of $F$, and let $G_F$ denote the Galois group $\mathrm{Gal}(F^s/F)$. Suppose $A$ is an abelian variety over $F$, and $\ell$ is a prime number not equal to the characteristic of $F$. Let $\bar{A} = A \times_F F^s$, let $A^t$ denote the dual of $A$, and let $\mathrm{End}^0(A) = \mathrm{End}(A) \otimes_{\mathbf{Z}} \mathbf{Q}$. Let $T_\ell(A)$ denote the $\ell$-adic Tate module $\varprojlim A_{\ell^m}$, where $A_{\ell^m}$ is the kernel of multiplication by $\ell^m$ in $A(F^s)$, and let $V_\ell(A) = T_\ell(A) \otimes_{\mathbf{Z}_\ell} \mathbf{Q}_\ell$. Let $\mathbf{Z}_\ell(1)$ denote the projective limit of the groups of $\ell^m$-th roots of unity in $\bar{\mathbf{F}}_p$.

To every invertible sheaf $\mathcal{L}$ on $\bar{A}$ corresponds a certain natural homomorphism of abelian varieties $\phi_{\mathcal{L}} : \bar{A} \to \bar{A}^t$, as defined in §13 of [14]. Let $\chi(\mathcal{L})$ denote the Euler characteristic of $\mathcal{L}$. Then $\deg(\phi_{\mathcal{L}}) = \chi(\mathcal{L})^2$ (see §16 of [14]), and $\phi_{\mathcal{L}}$ is an isogeny if and only if $\chi(\mathcal{L}) \neq 0$. If $\mathcal{L}$ is ample, then $\chi(\mathcal{L}) \neq 0$ and the isogeny $\phi_{\mathcal{L}}$ is called a polarization on $\bar{A}$. If the isogeny $\phi_{\mathcal{L}}$ is defined over $F$ (i.e., is obtained by extension of scalars from a morphism $A \to A^t$), then $\phi_{\mathcal{L}}$ is called an $F$-polarization on $A$. Now, assume that $\mathcal{L}$ is an invertible sheaf on $\bar{A}$ such that $\phi_{\mathcal{L}}$ is defined over $F$. Then $\mathcal{L}$ induces an alternating $G_F$-invariant bilinear form

$$E^{\mathcal{L}} : T_\ell(A) \times T_\ell(A) \to \mathbf{Z}_\ell(1) \cong \mathbf{Z}_\ell$$

as in §20 of [14]. The form $E^{\mathcal{L}}$ is nondegenerate if and only if $\phi_{\mathcal{L}}$ is an isogeny, i.e., $\chi(\mathcal{L}) \neq 0$. It is perfect if and only if $\phi_{\mathcal{L}}$ is an isogeny of degree prime to $\ell$. Thus, $E^{\mathcal{L}}$ is perfect if and only if $\chi(\mathcal{L})$ is not divisible by $\ell$.

Recall the finite group $N_p$ of order $2p(p-1)$ from Definition 2.8. Note that by Lemma 3.14, there exists a Galois extension $L$ of $\bar{\mathbf{F}}_p(t)$, totally ramified at $t = 0$, such that $\mathrm{Gal}(L/\bar{\mathbf{F}}_p(t)) \cong N_p \times \boldsymbol{\mu}_\ell$.

**Theorem 7.1.** *Suppose $p$ and $\ell$ are prime numbers, and $\ell$ does not divide $p(p-1)$. Let $r = (p-1)/2$ if $p$ is odd, and let $r = 1$ if $p = 2$. Suppose $b$ is a nonnegative integer, and $E$ is a supersingular elliptic curve over $\bar{\mathbf{F}}_p$. Let $F = \bar{\mathbf{F}}_p(t)$, let $F_0 = \bar{\mathbf{F}}_p((t))$, and let $L$ be a Galois extension of $F$, totally ramified at $t = 0$, such that $\mathrm{Gal}(L/F) \cong N_p \times \boldsymbol{\mu}_\ell$. Then there is an injective homomorphism*

$$c : \mathrm{Gal}(L/F) \hookrightarrow \mathrm{Aut}(E^{r(\ell-1)})$$

*such that if $A$ is the twist of $E^{r(\ell-1)}$ by $c$, and $B$ is an abelian variety over $F_0$ which is $F_0$-isogenous to $A \times E^b$, then $\chi(\mathcal{L})$ is divisible by $\ell$ for every invertible sheaf $\mathcal{L}$ on $\bar{B}$ with the property that $\phi_\mathcal{L}$ is defined over $F_0$. In particular, every $F_0$-polarization on $B$ has degree divisible by $\ell$.*

*Proof.* We begin by constructing an injective homomorphism $N_p \hookrightarrow \mathrm{Aut}(E^r)$, with respect to which $V_\ell(E^r)$ is an absolutely simple faithful symplectic $\mathbf{Q}_\ell[N_p]$-module.

If $p$ is odd, let $C$ be the hyperelliptic curve $y^2 = x^p - x$, and let $J(C)$ denote its Jacobian variety. Then $J(C) \cong E^r$ (see p. 172 of [6]). The subgroup of $\mathrm{Aut}(C)$ generated by $\zeta : (x,y) \mapsto (x+1, y)$ and

$$u : (x,y) \mapsto (\alpha^2 x, \alpha y), \quad \text{for some } \alpha \in \mathbf{F}_{p^2} \text{ with } \alpha^{p-1} = -1,$$

is isomorphic to $N_p$. Thus, $N_p \subseteq \mathrm{Aut}(C) \subseteq \mathrm{Aut}(J(C)) \cong \mathrm{Aut}(E^r)$. Since $\zeta^p = 1$ and $\zeta \neq 1$, the $\mathbf{Q}$-subalgebra $\mathbf{Q}(\zeta)$ of $\mathrm{End}^0(J(C))$ is isomorphic to either $\mathbf{Q}(\zeta_p)$ or $\mathbf{Q} \oplus \mathbf{Q}(\zeta_p)$. Since $\dim(J(C)) = (p-1)/2$, every commutative semisimple $\mathbf{Q}$-subalgebra of $\mathrm{End}^0(J(C))$ has dimension at most $p-1 = [\mathbf{Q}(\zeta_p) : \mathbf{Q}]$. Thus, $\mathbf{Q}(\zeta) \cong \mathbf{Q}(\zeta_p)$. Also, $u^{p-1} = -1 \in \mathrm{Aut}(J(C)) \cong \mathrm{Aut}(E^r)$, and $u$ acts (by conjugation) on $\mathbf{Q}(\zeta)$ as a generator of the Galois group $\mathrm{Gal}(\mathbf{Q}(\zeta)/\mathbf{Q})$. Therefore, the centralizer of $u$ in $\mathbf{Q}(\zeta) \otimes_{\mathbf{Q}} \mathbf{Q}_\ell$ is $\mathbf{Q}_\ell$. We have

$$\mathbf{Q}(\zeta) \otimes_{\mathbf{Q}} \mathbf{Q}_\ell \subset \mathrm{End}^0(J(C)) \otimes_{\mathbf{Q}} \mathbf{Q}_\ell \hookrightarrow \mathrm{End}_{\mathbf{Q}_\ell}(V_\ell(J(C))) \cong \mathrm{M}_{p-1}(\mathbf{Q}_\ell).$$

By dimension arguments, $\mathbf{Q}(\zeta) \otimes_{\mathbf{Q}} \mathbf{Q}_\ell$ is a maximal commutative subalgebra in $\mathrm{End}_{\mathbf{Q}_\ell}(V_\ell(J(C)))$. Since $N_p$ is generated by $\zeta$ and $u$, the centralizer of $N_p$ in $\mathrm{End}(V_\ell(J(C)))$ consists of scalars, i.e., the natural representation of $N_p$ on $V_\ell(J(C))$ ($\cong V_\ell(E^r)$) is absolutely irreducible. The theta divisor on $J(C)$ gives rise to an $N_p$-invariant principal polarization on $E^r$. We can therefore view $V_\ell(E^r)$ as an absolutely simple faithful symplectic $\mathbf{Q}_\ell[N_p]$-module.

If $p = 2$, then $N_2 \subset \mathrm{Aut}(E) \subset \mathrm{SL}(V_\ell(E))$ (see [23], §2). Since $N_2$ is a nonabelian finite group and $V_\ell(E)$ is 2-dimensional, therefore $V_\ell(E)$ is an absolutely simple faithful symplectic $\mathbf{Q}_\ell[N_2]$-module.

Let $G = N_p \times \boldsymbol{\mu}_\ell$. Fixing an isomorphism of additive groups $\mathbf{Z}[\zeta_\ell] \cong \mathbf{Z}^{\ell-1}$ gives rise to an embedding $\boldsymbol{\mu}_\ell \hookrightarrow \mathrm{GL}_{\ell-1}(\mathbf{Z})$. We thus have

$$\mathrm{Gal}(L/F) \cong G = N_p \times \boldsymbol{\mu}_\ell \subset \mathrm{Aut}(E^r) \times \mathrm{GL}_{\ell-1}(\mathbf{Z}) \subset \mathrm{Aut}(E^{r(\ell-1)}),$$

with $\mathrm{Aut}(E^r)$ embedded diagonally in $\mathrm{Aut}((E^r)^{\ell-1})$, and with the inclusion $\mathbf{Z} \subset \mathrm{End}(E^r)$ inducing the embedding $\mathrm{GL}_{\ell-1}(\mathbf{Z}) \subset \mathrm{Aut}((E^r)^{\ell-1})$. Denote the composition by

$$c : \mathrm{Gal}(L/F) \hookrightarrow \mathrm{Aut}(E^{r(\ell-1)}).$$

Let $A$ be the twist of $E^{r(\ell-1)}$ by the cocycle $c$. Let $L_0$ be the completion of $L$ at $t = 0$. Suppose $B$ is an abelian variety over $F_0$, and $\beta : B \to A \times E^b$ is an $F_0$-isogeny. By functoriality, $\beta$ induces an isomorphism

$$V_\ell(B) \cong V_\ell(A \times E^b) = V_\ell(A) \oplus V_\ell(E)^b$$

of $\mathbf{Q}_\ell[G_{F_0}]$-modules. Since $E$ is defined over the algebraically closed field $\bar{\mathbf{F}}_p$, therefore the $\ell$-power torsion on $E$ is defined over $\bar{\mathbf{F}}_p$, and thus over $F_0$. Therefore, $V_\ell(E^r)$ is a trivial $\mathbf{Q}_\ell[G_{F_0}]$-module. It follows that the $\ell$-adic representation $G_{F_0} \to \text{Aut}(T_\ell(A))$ factors through $\text{Gal}(L_0/F_0) \cong G$. ¿From the definitions of $c$ and $A$, it follows that

$$V_\ell(A) \cong V_\ell(E^r) \otimes_{\mathbf{Q}_\ell} \mathbf{Q}_\ell(\zeta_\ell)$$

as $\mathbf{Q}_\ell[G]$-modules (here we view $V_\ell(E^r)$ as an $N_p$-module and $\mathbf{Q}_\ell(\zeta_\ell)$ as a $\boldsymbol{\mu}_\ell$-module). Suppose $\mathcal{L}$ is an invertible sheaf on $\bar{B}$ such that $\phi_\mathcal{L}$ is defined over $F_0$ and $\chi(\mathcal{L})$ is not divisible by $\ell$. Then $E^\mathcal{L}$ is a perfect alternating $G$-invariant pairing on the lattice $T_\ell(B)$ in

$$V_\ell(B) \cong (V_\ell(E^r) \otimes_{\mathbf{Q}_\ell} \mathbf{Q}_\ell(\zeta_\ell)) \oplus V_\ell(E)^b.$$

This contradicts Theorem 6.2 (when $b = 0$) and Corollary 6.4, with $N = N_p$, $K = \mathbf{Q}_\ell$, $W = V_\ell(E^r)$, and $U = V_\ell(B)$. □

**Remark 7.2.** One may easily check that $L_0$ is the smallest extension of $F_0$ over which $B$ acquires good reduction [20].

**Corollary 7.3.** *Suppose $p$ and $\ell$ are prime numbers, $d$ is an integer, and $\ell$ does not divide $p(p-1)$. Let $r = (p-1)/2$ if $p$ is odd, and let $r = 1$ if $p = 2$. Suppose $d \geq r(\ell-1)$. Then there exists a $d$-dimensional supersingular abelian variety defined over $\bar{\mathbf{F}}_p(t)$ such that if $B$ is in its $\bar{\mathbf{F}}_p((t))$-isogeny class, then $\chi(\mathcal{L})$ is divisible by $\ell$ for every invertible sheaf $\mathcal{L}$ on $\bar{B}$ with the property that $\phi_\mathcal{L}$ is defined over $\bar{\mathbf{F}}_p((t))$. In particular, every $\bar{\mathbf{F}}_p((t))$-polarization on $B$ has degree divisible by $\ell$.*

**Remark 7.4.** Over every algebraically closed field of characteristic $p$, every supersingular abelian variety admits a polarization whose degree is a power of $p$ (see pp. 70–71 of [12]). Over finite fields, E. Howe showed that the isogeny class of every simple odd-dimensional abelian variety contains a principally polarized abelian variety (see Theorem 1.2 of [9]).

Department of Mathematics, Ohio State University, Columbus, Ohio 43210, USA
*E-mail address*: `silver@math.ohio-state.edu`

Department of Mathematics, Pennsylvania State University, University Park, PA 16802, USA
*E-mail address*: `zarhin@math.psu.edu`